\documentclass[journal]{IEEEtran}

\usepackage{subfigure}
\usepackage{cite}

\usepackage[cmex10]{amsmath}
\usepackage{graphicx}
\usepackage{amsfonts}

\newcommand{\mykeywords}[1]{\begin{IEEEkeywords} #1 \end{IEEEkeywords}}
\newcommand{\figsizea}[0]{6cm}
\newcommand{\figsizeb}[0]{5.5cm}

\begin{document}

\title{A compressive sensing approach for enhancing breast cancer detection using a hybrid DBT / NRI configuration}

\author{Richard Obermeier and Jose~Angel~Martinez-Lorenzo}

\maketitle



\begin{abstract}
This work presents a novel breast cancer imaging approach that uses compressive sensing in a 
hybrid Digital Breast Tomosynthesis (DBT) / Nearfield Radar Imaging (NRI) system configuration. The non-homogeneous tissue distribution of the breast, described in terms of  dielectric constant and conductivity, is extracted from the DBT image, and it is used by a full-wave Finite Difference in the Frequency Domain (FDFD) method to build a linearized model of the non-linear NRI imaging problem. The inversion of the linear problem is solved using compressive sensing imaging techniques, which lead to a reduction on the required number of sensing antennas and operational bandwidth without loss of performance. 

\mykeywords{breast cancer detection, microwave imaging, compressive sensing, convex optimization.}

\end{abstract}


\section{Introduction}
It has been shown \cite{Kopans2011} that X-ray based Digital Breast Tomosynthesis (DBT) enhances
the imaging capabilities of Conventional Mammography (CM) by producing volumetric imaging. 
Unfortunately, DBT still suffers from the limited radiological contrast existing between
cancerous and healthy fibroglandular tissue, which is of the order of 1\%. 
At microwave frequencies, the contrast between cancerous and healthy fibroglandular tissue is of
the order of 10\% \cite{Lazebnik2007}. For this reason, Nearfield Radar Imaging (NRI), 
working at these frequencies, is an appealing technology to detect cancerous tissue within the breast. 
Unfortunately, NRI by itself fails to detect cancerous tissue embedded in a random, heterogeneous
matrix of fibroglandular and fatty breast tissue.  


A recent paper \cite{MartinezLorenzo2013} demonstrated that it is possible to fuse DBT with NRI to achieve improved breast cancer detection. 
The successful detection of breast cancer using this configuration required a
high performance NRI sensor operating with the following requirements: 1) $17$ transmitting/receiving antennas in a multi-static configuration, and 2) a wideband system, having a $1$GHz bandwidth at a $1.5$GHz center frequency. These high performance requirements make the NRI sensor too expensive for widespread use in a clinical setting.

This paper presents preliminary 2D results of a new 3D imaging approach for
the NRI sensor \cite{molaei2015}, based on techniques from compressive 
sensing theory. This approach utilizes a multiple monostatic configuration with fewer antennas and a smaller operating bandwidth than the NRI system presented in \cite{MartinezLorenzo2013}.
The new approach reduces the cost of the NRI system without loss of performance, thereby making the technology more suitable for clinical use.

\section{Methodology}
\label{sec:method}
\subsection{System Configuration}
The concept of a fused DBT / NRI imaging system was first
introduced in \cite{MartinezLorenzo2013, Obermeier2014}, and the configuration of the 3D imaging mechatronics device was presented in \cite{molaei2015}. In this approach,
the breast is placed under clinical compression, and the DBT measurements are recorded using
low-dosage X-rays. The DBT measurements are used to create a 3D
reconstruction of the fat distribution in the breast. Simultaneously, a set of NRI
measurements are collected using a series of microwave antennas,
which are placed along the breast in a bolus material in order to maximize its coupling with the breast surface. 

\subsection{The Sensing Problem}
\label{sec:sensing_form}
The NRI system operates using a multiple monostatic configuration, in which each antenna utilizes stepped-frequency waveforms. Measurements of $n_f$ frequencies are
taken by $n_a$ antennas for a total of $M =n_a\cdot n_f$ measurements. 
The relationship between the electric fields $\mathbf{E}(\mathbf{r},\omega)$ and the unknown complex permittivity $\epsilon(\mathbf{r}, \omega)$ of the breast tissues can be expressed as:
\begin{align} 
\mathbf{E}_s(\mathbf{r},\omega) &=  \int \mathbf{G}_b(\mathbf{r}, \mathbf{r}', \omega) k_b^2(\mathbf{r}', \omega)\mathbf{E}(\mathbf{r}',\omega)\chi(\mathbf{r}',\omega) d\mathbf{r}' \label{eq:cont_src} \\
\mathbf{E}(\mathbf{r},\omega) &= \mathbf{E}_b(\mathbf{r},\omega) + \mathbf{E}_s(\mathbf{r}, \omega)
\end{align}
where $\mathbf{G}_b(\mathbf{r}, \mathbf{r}', \omega)$ are the Green's functions of the background medium, $\mathbf{E}_b(\mathbf{r},\omega)$ are the incident electric fields, $\mathbf{E}_s(\mathbf{r},\omega)$ are the scattered electric fields, and $\chi(\mathbf{r},\omega)=\frac{\epsilon(\mathbf{r},\omega)-\epsilon_b(\mathbf{r},\omega)}{\epsilon_b(\mathbf{r},\omega)}$ are the contrast variables \cite{VanDenBerg2001,Abubakar2011,Zakaria2013}.

Eq. \ref{eq:cont_src} is a nonlinear function of the contrast variable $\chi(\mathbf{r}, \omega)$
and total electric field $\mathbf{E}(\mathbf{r}, \omega)$, and so nonlinear programming techniques such as the Contrast Source Inversion (CSI) algorithm \cite{VanDenBerg2001,Abubakar2011,Zakaria2013} must be applied in order to recover $\chi(\mathbf{r}, \omega)$. These types of nonlinear algorithms typically require several calls to a forward model solver in each iteration, which makes them computationally expensive. As a result, it is desirable to make some simplifying assumptions in order to reduce the computation time. This work makes two such assumptions. First, the Born approximation is applied, $\mathbf{E}(\mathbf{r},\omega) \approx \mathbf{E}_b(\mathbf{r},\omega)$, in order to linearize Eq. \ref{eq:cont_src}. Second, the complex permittivities are assumed to be approximately constant over the frequency range of the NRI system, i.e. $\epsilon(\mathbf{r}, \omega) \approx \epsilon(\mathbf{r})$ and $\epsilon_b(\mathbf{r}, \omega) \approx \epsilon_b(\mathbf{r})$, so that the contrast variable is also approximately constant over frequency. With these two modifications, Eq. \ref{eq:cont_src} can be rewritten in the following form:
\begin{align}
\mathbf{E}_s(\mathbf{r},\omega) &=  \int\mathbf{G}_b(\mathbf{r}, \mathbf{r}', \omega) k_b^2(\mathbf{r}', \omega)\mathbf{E}_b(\mathbf{r}',\omega)\chi(\mathbf{r}') d\mathbf{r}'  \label{eq:born_approx} \\
&+ \mathbf{\hat{e}}_s(\mathbf{r}, \omega)\nonumber
\end{align}
where $\mathbf{\hat{e}}_s(\mathbf{r}, \omega)$ is the error introduced by the approximating assumptions. 

Eq. \ref{eq:born_approx} can be discretized as $\mathbf{y} = \mathbf{A}\mathbf{x} + \mathbf{e} + \boldsymbol{\eta}$, where $\mathbf{x} \in \mathbb{C}^N$ are the contrast variables, $\mathbf{y} \in \mathbb{C}^M$ are the measured fields, $\mathbf{A} \in \mathbb{C}^{M\times N}$ is the sensing matrix constructed from the incident fields and Green's functions of the background medium, $\mathbf{e} \in \mathbb{C}^M$ is the error vector, and $\boldsymbol{\eta} \in \mathbb{C}^M$ is the random noise introduced by the measurement system. In practice, $M < N$, and so this system has an infinite number of solutions satisfying $\mathbf{y} = \mathbf{A}\mathbf{x} + \mathbf{e} + \boldsymbol{\eta}$. When $\|\mathbf{e}\|_{\ell_2} \ll \|\boldsymbol{\eta}\|_{\ell_2}$, the performance of linear inverse techniques only depends upon the vector $\boldsymbol{\eta}$. When $\|\mathbf{e}\|_{\ell_2} \sim \|\boldsymbol{\eta}\|_{\ell_2}$, then the performance of linear inverse techniques depends upon both $\mathbf{e}$ and $\boldsymbol{\eta}$. In practice, the statistics of the measurement noise $\boldsymbol{\eta}$ can be estimated in order to tune the proposed inverse algorithm accordingly. However, it is difficult to estimate $\mathbf{e}$, since it requires a-priori knowledge of the unknown contrast variable $\mathbf{x}$.

Stand-alone NRI systems tend to perform poorly in breast cancer imaging applications because they lack a suitable background model $\epsilon_b(\mathbf{r},\omega)$. Without any prior knowledge, NRI systems select a homogeneous background medium that has a dielectric constant derived from averaging that of low-water-content (LWC) fatty tissue and high-water-content (HWC) fibroglandular tissue. Choosing a homogeneous dielectric constant leads to a contrast variable $\mathbf{x}$ that has a significant number of large, non-zero elements. This violates the assumptions made by the Born approximation and produces an error vector $\mathbf{e}$ with a large norm, which ultimately challenges the accurate inversion of the linear system of equations.

\begin{figure}[h!]
     \begin{center}
            \includegraphics[height=\figsizea]{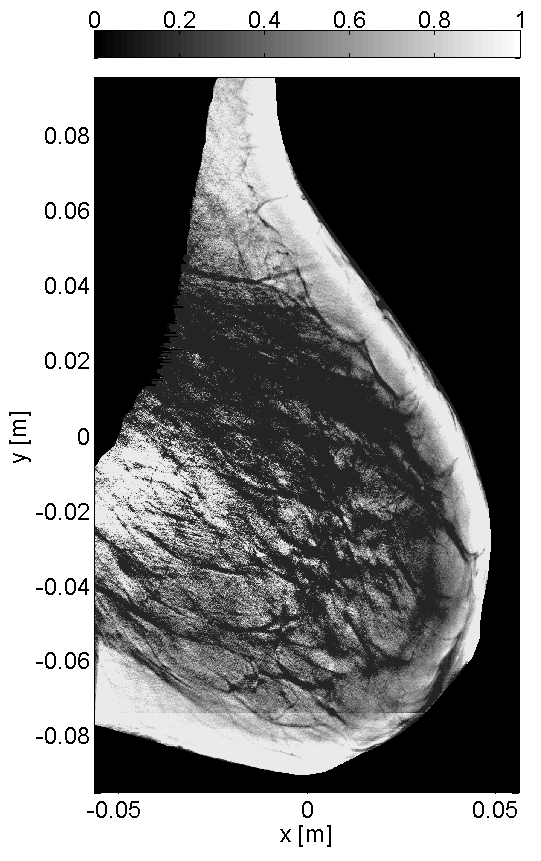}
    \caption{Fat content of breast segmented from a 3D DBT image. High intensity indicates a high percentage of fat.}
   \label{fig:segment}
\end{center}
\end{figure}

The hybrid DBT / NRI system utilizes the prior knowledge obtained from the DBT image in order to construct the heterogeneous background model, thereby overcoming the aforementioned limitations of the homogeneous model used by stand-alone NRI systems. In this hybrid system, the DBT image is segmented into three types of tissue, skin, muscle (pectoralis major), and breast tissue, and it is assumed that the latter only contains healthy tissue. Additionally, each pixel of breast tissue consists of $p\%$ of fatty tissue and $(100-p)\%$  of fibroglandular tissue, and its constitutive properties $\epsilon_r(\mathbf{r}, \omega)$ and $\sigma(\mathbf{r}, \omega)$ are determined from a composite model, which considers the dispersive properties of the breast tissues \cite{MartinezLorenzo2013a}.
Fig. \ref{fig:segment} displays the fat content segmented from a single slice of a 3D DBT image. The segmented geometry is input to an electromagnetic numerical simulation based on Finite Differences in the Frequency Domain (FDFD) \cite{Rappaport2001} in order to model the NRI sensing process. The FDFD is a full-wave model that accounts for all mutual interactions within the breast, and is used to generate the incident fields $\mathbf{E}_b(\mathbf{r}, \omega)$ and the non-uniform, cancer free Green's functions $\mathbf{G}_b(\mathbf{r}, \mathbf{r}', \omega)$.

\section{Imaging with Compressive Sensing}
Ideally, the hybrid DBT / NRI system would segment the healthy breast tissue perfectly and would classify any cancerous tissue as HWC fibroglandular tissue, producing a contrast variable $\mathbf{x}$ that is non-zero only at the locations of cancerous lesions. In this case, the reconstruction process becomes a sparse recovery problem, which can be solved using techniques from compressive sensing (CS) theory. CS theory is a relatively novel signal processing technique, which was first introduced by Candes et al. in \cite{Candes2006a}, and it has since been refined in many other works \cite{Candes2006,Donoho2006,Becker2011,massa2015}. CS establishes that sparse signals can be recovered using far fewer measurements than required by the Nyquist sampling criterion. In order to apply such principles, certain mathematical conditions must be satisfied by the sensing matrix $\mathbf{A}$ and the reconstructed image $\mathbf{x}$. The sensing matrix must satisfy the Restricted Isometry Property condition, which is related to the independence of its columns, and the number of non-zero entries, $N_{nz}$, of the reconstructed image must be much smaller than the total number of elements, $N$. If the two aforementioned conditions are satisfied, then the reconstruction of the unknown vector can be performed with a small number of measurements by solving a norm-$1$ optimization problem. In this work, the contrast variables $\mathbf{x}$ are recovered by solving the modified basis pursuit denoising problem:
\begin{align}
\underset{\mathbf{x}}{\text{minimize}} ~~&\frac{1}{2}\|\mathbf{A}\mathbf{x}-\mathbf{y}\|_{\ell_2}^2 + \lambda \|\mathbf{x}\|_{\ell_1} \label{eq:l1opt} \\
\text{subject to} ~~&\operatorname{\mathbf{Re}}(\operatorname{diag}({\boldsymbol\epsilon}_b)\mathbf{x} + {\boldsymbol\epsilon}_b) \ge \mathbf{1} \nonumber  \\
~~&\operatorname{\mathbf{Im}}(\operatorname{diag}({\boldsymbol\epsilon}_b)\mathbf{x} + {\boldsymbol\epsilon}_b) \ge \mathbf{0}  \nonumber
\end{align}
where the constraints on the variable $\mathbf{x}$ ensure that the solution is physically realizable, i.e. $\epsilon_r \ge 1$ and $\sigma \ge 0$. In practice, the weighting factor $\lambda$ must be selected based upon the expected error in the measurement vector $\mathbf{y}$ and sparsity of the vector $\mathbf{x}$. This problem can be efficiently solved using Nesterov's accelerated gradient method for non-smooth convex functions \cite{Nesterov2005}, which was first applied to the basis pursuit denoising problem in \cite{Becker2011}. In this method, the norm-$1$ term is replaced by the smooth approximation:
\begin{equation}
g_\mu(\mathbf{x}) = g_\mu(x_1, \ldots, x_N) = \sum_{n=1}^N f_\mu(x_n)
\end{equation}
where the Huber function $f_\mu(\cdot)$ is defined as:
\begin{equation}
f_\mu(x) = \begin{cases}
|x| & |x| \ge \mu \\
\frac{1}{2\mu}|x|^2 & |x| < \mu
\end{cases}
\end{equation}
The Huber function is Lipschitz continuously differentiable, so first-order techniques such as Nesterov's method can be used to minimize the smoothed objective function $\frac{1}{2}\|\mathbf{A}\mathbf{x}-\mathbf{y}\|_{\ell_2}^2 + \lambda g_\mu(\mathbf{x})$. At each step in Nesterov's method, the desired variable is updated using a step direction derived from the gradients up to and including the current iteration. This updated point is then projected onto the feasible set to ensure that all of the constraints are satisfied. For the imaging problem of Eq. \ref{eq:l1opt}, the projection operator can be expressed as the solution to the following convex optimization problem:
\begin{align}
\underset{\mathbf{x}}{\text{minimize}} ~~&\|\mathbf{x}-\mathbf{z}\|_{\ell_2}^2 \label{eq:prox} \\
\text{subject to} ~~&\operatorname{\mathbf{Re}}(\operatorname{diag}({\boldsymbol\epsilon}_b)\mathbf{x} + {\boldsymbol\epsilon}_b) \ge \mathbf{1} \nonumber  \\
~~&\operatorname{\mathbf{Im}}(\operatorname{diag}({\boldsymbol\epsilon}_b)\mathbf{x} + {\boldsymbol\epsilon}_b) \ge \mathbf{0}  \nonumber
\end{align}
This problem is separable in the elements $x_1,\ldots,x_N$ of the vector $\mathbf{x}$. For a scalar contrast $x$, it can be shown that this projection problem has the following closed form solution:
\begin{align}
P_{Q_p}(x) &= \frac{P_\epsilon(\epsilon_b x + \epsilon_b) - \epsilon_b}{\epsilon_b}\label{eq:proj1} \\
P_{\epsilon}(\epsilon) &= \operatorname{max}(\operatorname{Re}(\epsilon),1) + \jmath \operatorname{max}(\operatorname{Im}(\epsilon), 0) \label{eq:proj2}
\end{align}
The reader is referred to the literature \cite{Nesterov2005,Becker2011} for further details on Nesterov's method
\section{Numerical Results}
\label{sec:results}
Following the process of Section \ref{sec:sensing_form}, a 2D model of a healthy breast was generated by segmenting a 2D slice from a 3D  DBT image. In order to simulate data from a cancerous case, a lesion with frequency-dependent electrical properties modeled after \cite{Lazebnik2007} was added to the healthy breast. A 2D version of the FDFD code was used to generate the synthetic NRI measurements of the healthy breast, the synthetic NRI measurements of the cancerous breast, and the sensing matrix of the healthy breast $\mathbf{A}$ according to Eq. \ref{eq:born_approx}. Note that the FDFD model accounted for the dispersive properties of both the healthy breast tissue and the cancerous tissue; only the inversion process utilized the simplifying assumptions of Section \ref{sec:sensing_form}. In the simulation, the NRI system used six transmitting and receiving antennas operating in a multiiple monostatic configuration. Each antenna was excited with three different frequencies, $500$MHz,  $600$MHz, and $700$MHz, for a total of $18$ measurements among the antennas.
\begin{figure}[h!]
\centering
\includegraphics[height=\figsizeb, clip=true]{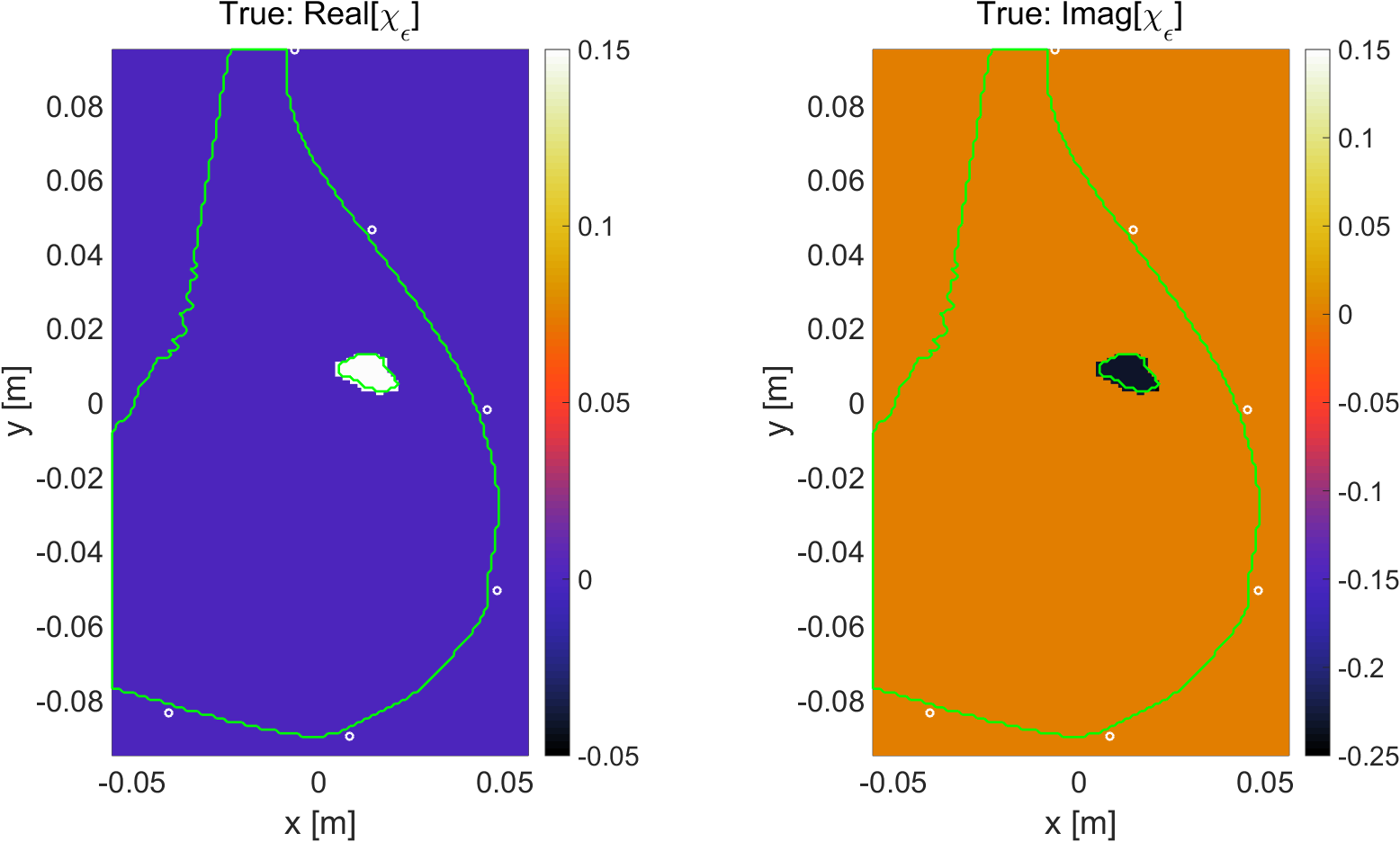}
    \caption{Real and imaginary parts of true contrast variable $\chi_{{\epsilon}}$ obtained when the DBT image is segmented perfectly.}
   \label{fig:true_contrast_000}
\end{figure}

\begin{figure}[h!]
\centering
\includegraphics[height=\figsizeb, clip=true]{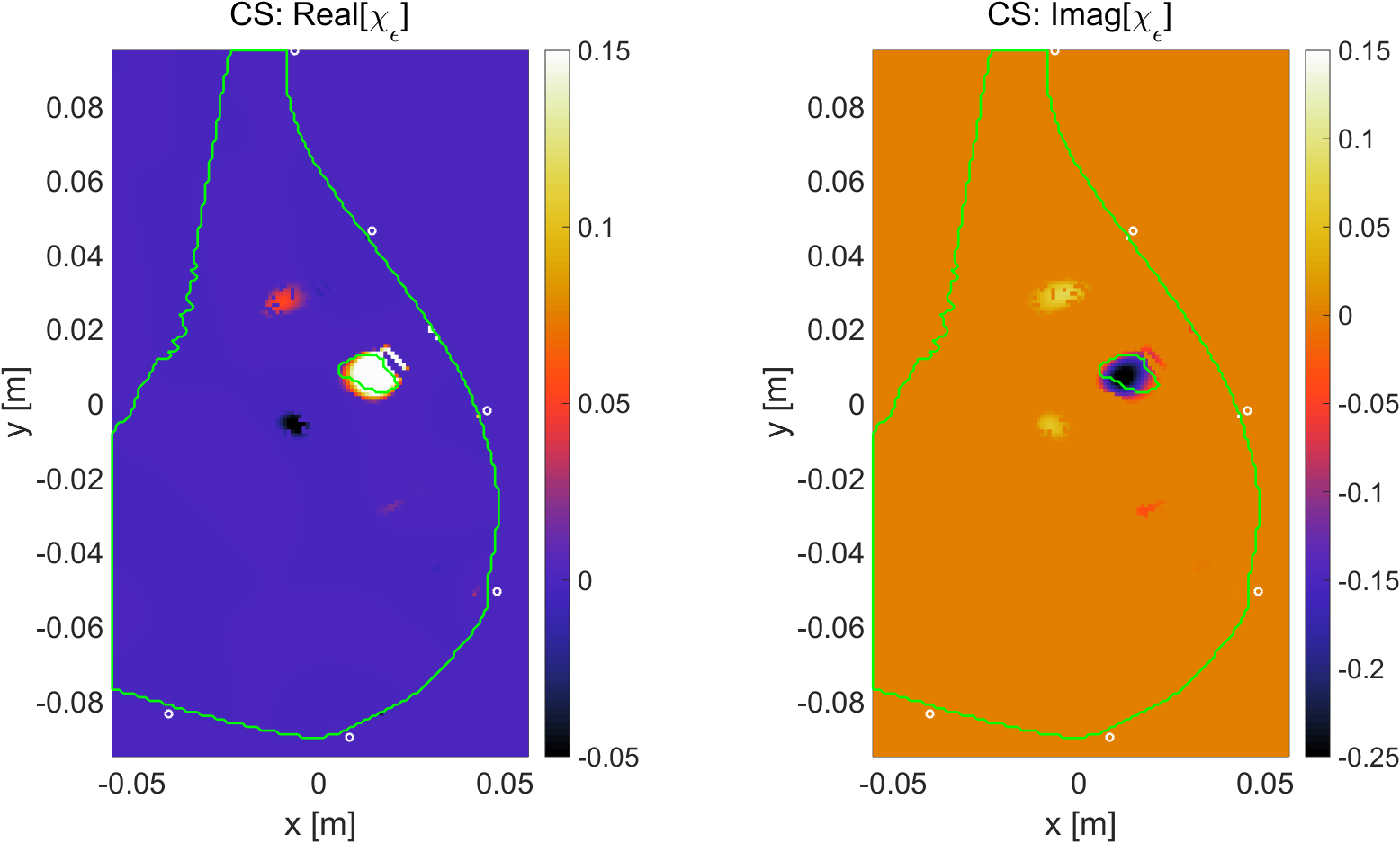}
    \caption{Real and imaginary parts of reconstructed contrast variable $\hat{\chi}_{{\epsilon}}$ obtained when the DBT image is segmented perfectly and there is no measurement noise.}
   \label{fig:cs_contrast_000_000}
\end{figure}

Figure \ref{fig:true_contrast_000} displays the true contrast variable obtained when the fat percentage is perfectly segmented from the DBT image. In this plot, the white dots represent the antenna positions and the green curves represent the breast and lesion borders. Since the fat percentage was segmented perfectly, the contrast variable is non-zero only at the location of the cancerous lesion. Figure \ref{fig:cs_contrast_000_000} displays the estimated contrast variable obtained using the perfect fat percentage segmentation and noiseless measurements. The artifacts within the image are due to the error vector $\mathbf{\hat{e}}_s(\mathbf{r}, \omega)$ that is introduced to the measurement vector when the simplifying assumptions of Section \ref{sec:sensing_form} are applied. Despite these artifacts, the algorithm is able to locate the cancerous lesion. This represents a significant improvement over the phase-based results presented in \cite{MartinezLorenzo2013}, which utilized $17$ antennas in a multi-static configuration and $11$ frequencies over a $1$GHz bandwidth in order image the cancerous lesion.



\begin{figure}[b!]
\centering
\includegraphics[height=\figsizeb, clip=true]{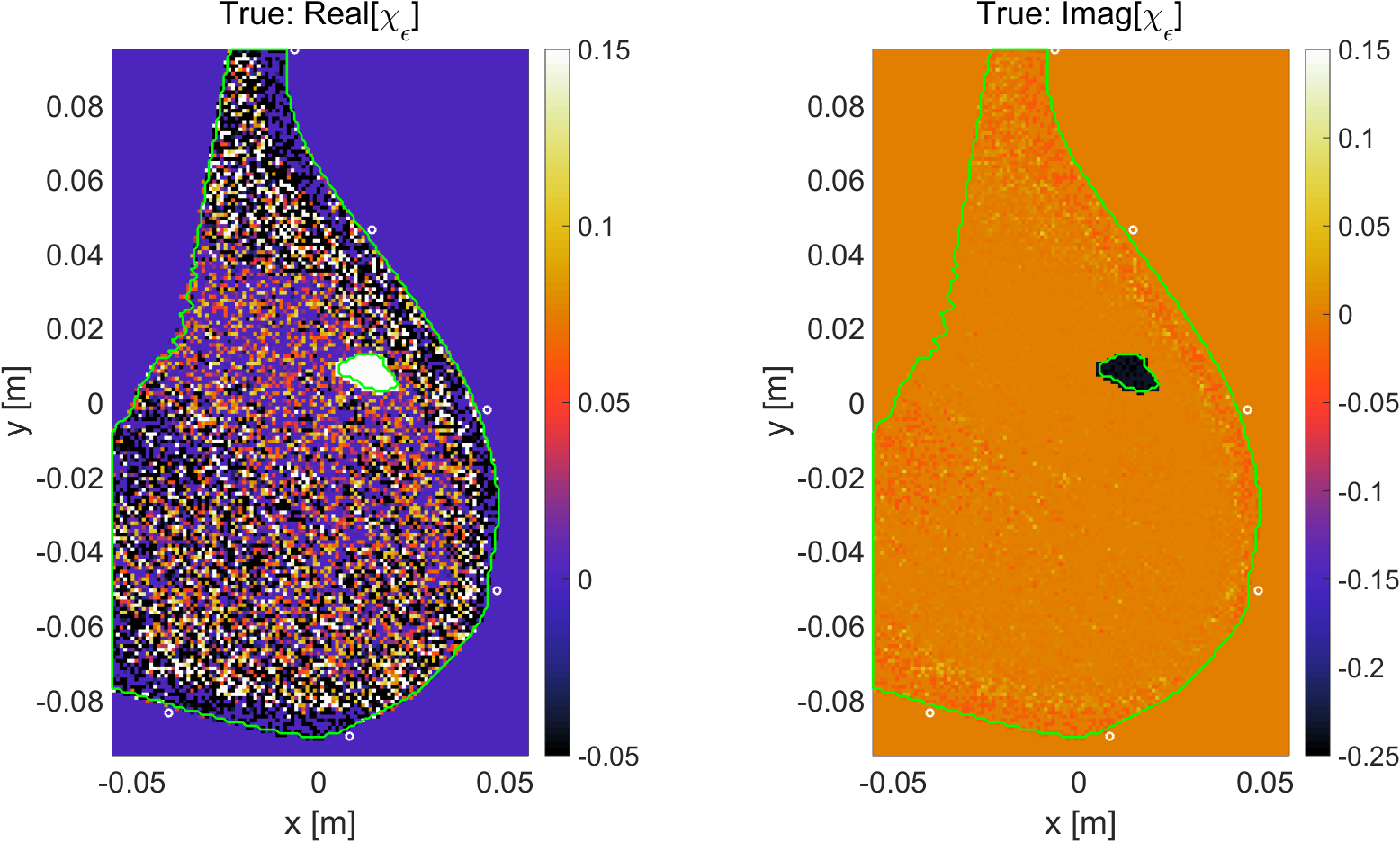}
    \caption{Real and imaginary parts of true contrast variable $\chi_{{\epsilon}}$ obtained when the fat percentage is segmented from the DBT image with $10\%$ error.}
   \label{fig:true_contrast_100}
\end{figure}

\begin{figure}[b!]
\centering
\includegraphics[height=\figsizeb, clip=true]{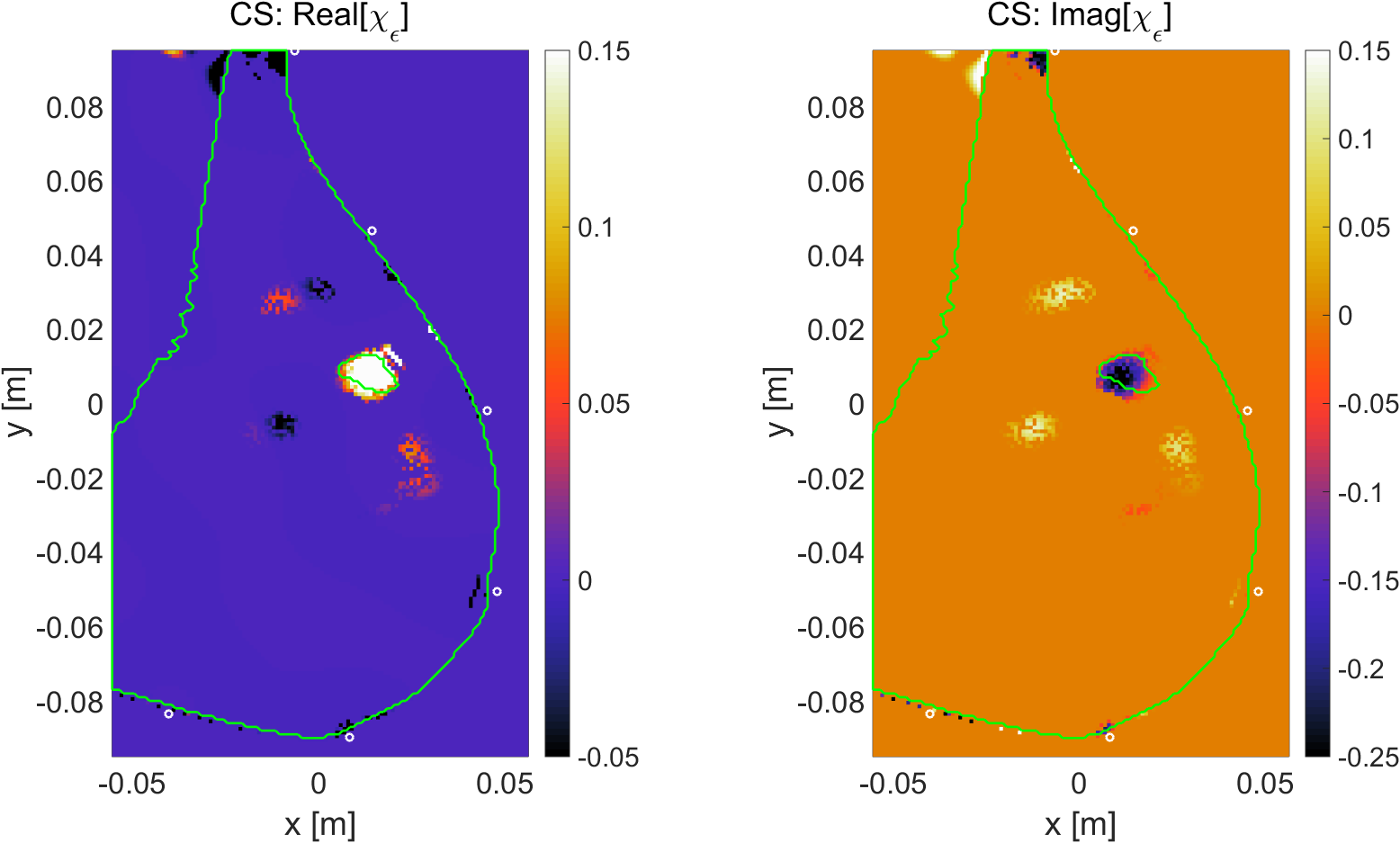}
    \caption{Real and imaginary parts of reconstructed contrast variable $\hat{\chi}_{{\epsilon}}$ obtained when the fat percentage is segmented from the DBT image with $10\%$ error and there is no measurement noise.}
   \label{fig:cs_contrast_100_000}
\end{figure}

Figure \ref{fig:true_contrast_100} displays the true contrast variable obtained when the fat percentage is segmented from the DBT image with $10\%$ random error. More specifically, the fat percentage values were corrupted by i.i.d. random noise following a uniform distribution, taking values between $\pm 10\%$ with equal probability. Since the fat percentage is not segmented correctly, the true contrast variable is non-zero within the healthy tissue. Nevertheless, the true contrast variable is approximately compressible, and so Eq. \ref{eq:l1opt} can still be used to image the breast. This result can be seen in Figure \ref{fig:cs_contrast_100_000}, which displays the estimated contrast variable obtained using the noisy fat percentage segmentation and noiseless measurements. 

\begin{figure}[h!]
\centering
\includegraphics[height=\figsizeb, clip=true]{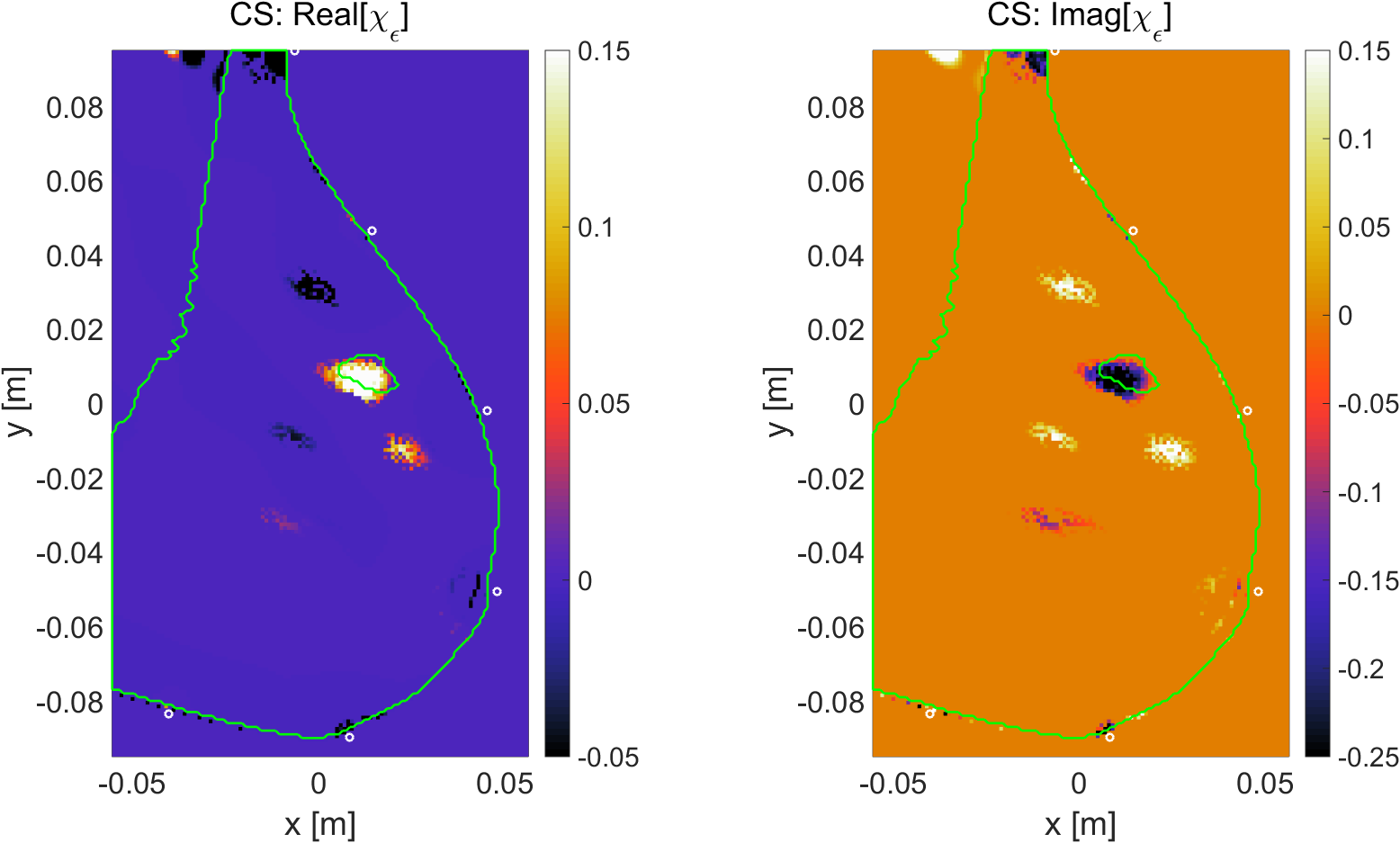}
    \caption{Real and imaginary parts of reconstructed contrast variable $\hat{\chi}_{{\epsilon}}$ obtained when the fat percentage is segmented from the DBT image with $10\%$ error and and the measurement SNR $=49\text{dB}$.}
   \label{fig:cs_contrast_100_300}
\end{figure}

\begin{figure}[h!]
\centering
\includegraphics[height=\figsizeb, clip=true]{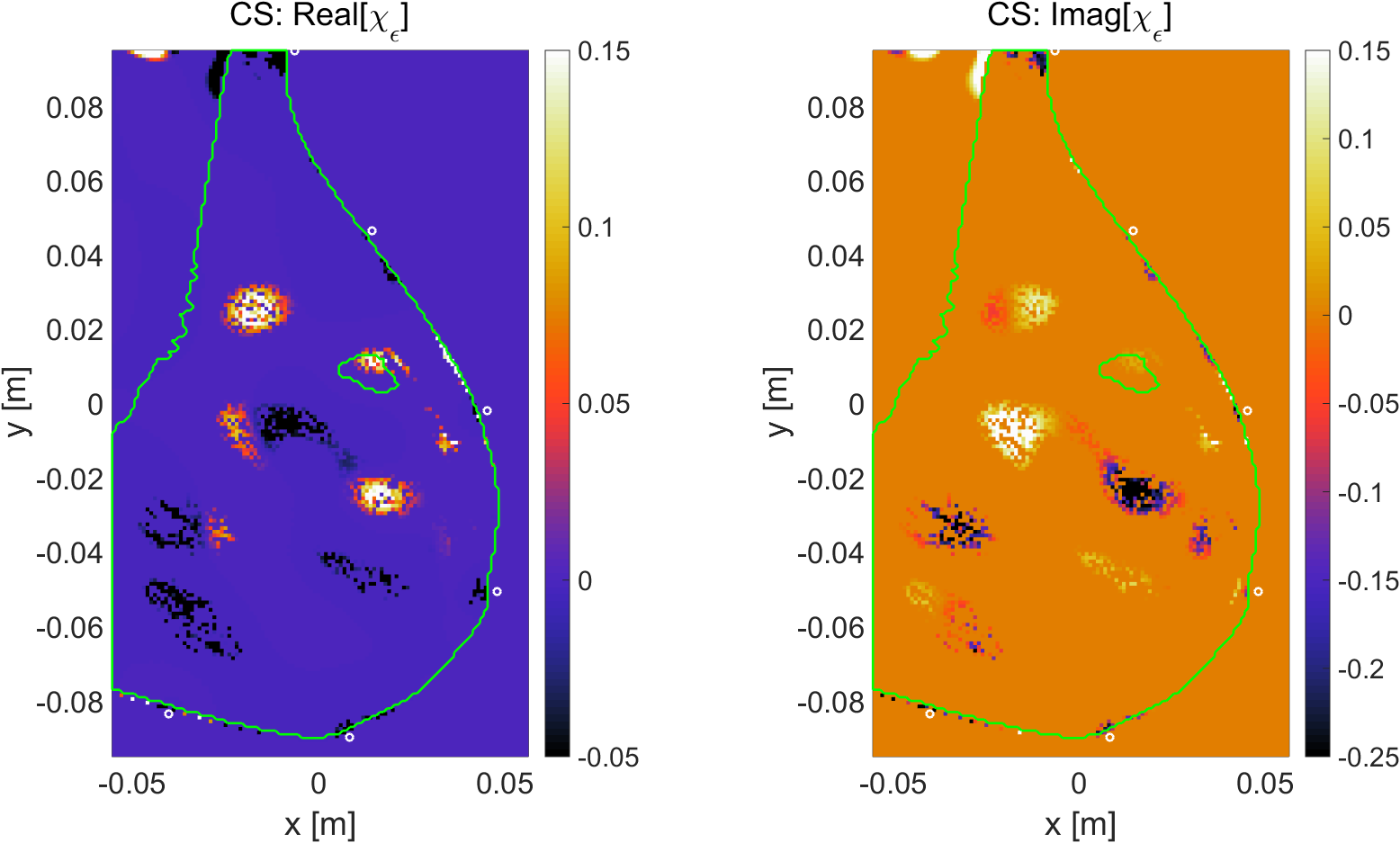}
    \caption{Real and imaginary parts of reconstructed contrast variable $\hat{\chi}_{{\epsilon}}$ obtained when the fat percentage is segmented from the DBT image with $10\%$ error and and the measurement SNR $=43\text{dB}$.}
   \label{fig:cs_contrast_100_500}
\end{figure}

Figure \ref{fig:cs_contrast_100_300} displays the estimated contrast variable obtained using the noisy fat percentage segmentation and measurements whose SNR $=49\text{dB}$. It is important to note that the signals in this SNR calculation are the electric fields scattered by the entire breast, and not just the fields scattered by the cancerous lesion. In this example, the fields scattered by the lesion are approximately $40\text{dB}$ lower in magnitude than the fields scattered by the rest of the breast, so that the ``lesion signal to noise ratio'' is on the order of $10\text{dB}$. Therefore, the NRI system must have a significant SNR to ensure that the fields produced by cancerous lesions are not overwhelmed by the noise, or it must have antennas with higher directivity in order to improve the SNR - the latter case may require the CS algorithm to use additional measurements. With this high SNR, the CS algorithm is able to image the cancerous lesion with some additional artifacts compared to the noiseless case. However, when the SNR is decreased to $43\text{dB}$, the algorithm is no longer able to image the lesion, as can be seen in Figure \ref{fig:cs_contrast_100_500}. 

\section{Conclusions}
\label{sec:conclusions}
This work presents a novel approach for imaging breast cancer using compressive sensing techniques in a hybrid DBT / NRI imaging system. Using the prior knowledge obtained from the DBT reconstruction, the hybrid system creates a background tissue distribution of the assumed healthy breast, thereby overcoming a common pitfall in stand-alone NRI imaging systems. Since cancerous lesions tend to be localized to a relatively small region of the breast, image reconstruction can be expressed as a sparse recovery problem, which is solved using techniques from compressive sensing. Numerical results show that the CS imaging algorithm can localize cancerous lesions within in the breast, even when corrupted by DBT segmentation and measurement errors.

\section*{Acknowledgment}
This work has been funded by the U.S. National Science Foundation award number 1347454.

\ifCLASSOPTIONcaptionsoff
  \newpage
\fi

\bibliography{./SICA-TA}
\bibliographystyle{IEEEtran}

\end{document}